# SPACE-TIME MULTI-PATCH DISCONTINUOUS GALERKIN ISOGEOMETRIC ANALYSIS FOR PARABOLIC EVOLUTION PROBLEMS

Stephen Edward Moore

Katholische Hochschulgemeinde der Diözese Linz,[**]
Petrinumstrasse 12/8, A-4040, Linz.

**Abstract.** We present and analyze a stable space-time multi-patch discontinuous Galerkin Isogeometric Analysis (dG-IgA) scheme for the numerical solution of parabolic evolution equations in deforming space-time computational domains. Following [20], we use a time-upwind test function and apply multi-patch discontinuous Galerkin IgA methodology for discretizing the evolution problem both in space and in time. This yields a discrete bilinear form which is elliptic on the IgA space with respect to a space-time dG norm. This property together with a corresponding boundedness property, consistency and approximation results for the IgA spaces yields an *a priori discretization* error estimate with respect to the space-time dG norm. The theoretical results are confirmed by several numerical experiments with low- and high-order IgA spaces.

## 1 Introduction

Several discretization methods exists for solving parabolic evolution problems. The standard discretization methods in time and space are based on time-stepping methods combined with some spatial discretization technique like the finite element methods, finite difference methods, finite volume methods,etc. The two main time-stepping methods are horizontal and vertical method of lines. The vertical method of lines discretize first in space and then in time, whereas in the horizontal method of lines, also known as Rothe's method, the discretization starts with respect to the time variable, see e.g. [17, 15, 28]. However, due to the separation of space and time discretizations, the development of efficient and fully adaptive schemes becomes complicated. This also affects negatively the parallelization of solvers due to the curse of sequentiality.

---

[**] Corresponding address: moorekwesi@gmail.com
This work was started when the author was working at RICAM, Altenberger Strasse 69, A-4040 Linz, Austria.

Another approach to discretization of parabolic evolution problems are the space-time methods. The space-time schemes allow for discretization in space and time simultaneously see [12, 30]. Space-time and also space-time discontinuous Galerkin (dG) discretization and their applications to several partial differential equations including their adaptivity and the treatment of changing spatial domains in time have been addressed for example, in [27, 5, 22] and [23] and all references therein.

Many approximation schemes have been developed for space-time formulation of parabolic evolution problems. Particularly, the $p$ and $hp$ in time versions of $hp-$finite element method to parabolic problems have been investigated in [2] and [3] respectively. This was followed by $hp$ dG time-stepping combined with FEM in space see [24]. Also, space-time wavelet methods for parabolic evolution problems have been analyzed in [25, 7]. Error bounds for reduced basis approximation see e.g. [29] and uniform stability of an abstract Petrov-Galerkin discretization of boundedly invertible operators for space-time discretizations to linear parabolic problems have been studied in [1]. Space-time streamline diffusion methods which uses special time-upwind test functions have also been applied to Navier-Stokes and other engineering problems see e.g. [10, 13].

Another numerical scheme that has been introduced in recent times is the Isogeometric Analysis (IgA). Isogeometric analysis was introduced by Hughes et al. in 2005, see, e.g. [11] and uses the same basis functions for both representing the approximate solution and the geometry description. In [4], the authors introduced approximation, stability and error estimates for $h$-refined IgA meshes of spatial computational domains see also [26]. The tools involved in the analysis and experimentation of the scheme are computational domain $Q$, a geometrical mapping $\Phi$ and a parameter domain $\widehat{Q}$. The computational domain also called patch can involve one or several domains. In most industrial applications, the computational domains involve several domains or patches. Several assembling strategies involving multiple patches have been presented in IgA including mortar methods [6] and discontinuous Galerkin (dG) methods [18].

In [18], the dG-IgA was presented where the approximation estimates are valid for each patch including non-matching meshes. Using these approximation results for B-splines or NURBS together with the elliptic-ity of the discrete bilinear form $a_h(\cdot, \cdot)$ with respect to a discrete norm $\| \cdot \|_h$ with boundedness and consistency results, asymptotically optimal discretization error estimate in the discrete dG-IgA norm $\| \cdot \|_h$ were obtained.

The aim of this paper is to present a new stable space-time multipatch dG-IgA scheme for the numerical approximation of parabolic initial boundary value problems. The current article is motivated by an earlier article, see [20], where a single-patch moving and also nonmoving spatial computational domain were considered. The space-time cylinder $Q$ was represented by one smooth, uniformly regular B-spline or NURBS map $\Phi$ of the parameter domain $\widehat{Q} = (0,1)^{d+1}, d \in \{1, 2, 3\}$. By using a time-upwind test function, we derive a discrete bilinear form $a_h(\cdot, \cdot) : V_{0h} \times V_{0h} \to \mathbb{R}$ which is elliptic on the space-time dG-IgA space $V_{0h}$ with respect to a discrete space-time norm $\| \cdot \|_h$.

The rest of the paper is organized as follows; Section 2, we introduce the model problem and variational formulation. We also present tools necessary for the IgA. In Section 3, we will present the derivation of the stable space-time multi-patch discontinuous Galerkin Isogeometric analysis. Section 4 will deal with the approximation results and finally, the *a priori error* estimate. We will present numerical results showing to confirm theory presented in Section 5

## 2   Model Problem and Method

Let us first introduce the Sobolev spaces $H^{s_x, s_t}(Q) = \{u \in L_2(Q) : \partial_x^\alpha u \in L_2(Q), \forall \alpha \text{ with } 0 \leq |\alpha| \leq s_x, \partial_t^i u \in L_2(Q), i = 0, \ldots, s_t\}$ of functions defined in the space-time cylinder $Q$, where $L_2(Q)$ denotes the space of square-integrable functions, $\alpha = (\alpha_1, ..., \alpha_d)$ is a multi-index with non-negative integers $\alpha_1, ..., \alpha_d$, $|\alpha| = \alpha_1 + \ldots + \alpha_d$, $\partial_x^\alpha u := \partial^{|\alpha|} u / \partial x^\alpha = \partial^{|\alpha|} u / \partial x_1^{\alpha_1} \ldots \partial x_d^{\alpha_d}$ and $\partial_t^i u := \partial^i u / \partial t^i$, see, e.g., [16].

### 2.1   Model Problem

We consider a linear parabolic initial-boundary value problem: find $u : \overline{Q} \to \mathbb{R}$ such that

$$\partial_t u - \Delta u = f \quad \text{in} \quad Q, \quad u = 0 \quad \text{on} \quad \Sigma \cup \Sigma_0, \qquad (2.1)$$

as the model problem posed in the space-time cylinder $\overline{Q} = \overline{\Omega(t)} \times [0, T]$, where $\Delta$ is the Laplace operator, $\partial_t$ denotes the partial time derivative, $f$ is a given source function, $u_0$ is the given initial data, $T$ is the final time, and $\Omega(t) \subset \mathbb{R}^d$, $d \in \{1, 2, 3\}$ denotes the spatial computational domain with the boundary $\partial\Omega$. We consider that the computational domain $Q$ is bounded and Lipschitz. The standard weak space-time variational formulation of (2.1) reads as follows: find $u \in H_{0,\underline{0}}^{1,0}(Q)$ such that

$$a(u, v) = l(v) \quad \forall v \in H_{0,\overline{0}}^{1,1}(Q), \tag{2.2}$$

with the bilinear and linear forms given by

$$a(u, v) = -\int_Q u\partial_t v dx dt + \int_Q \nabla_x u \cdot \nabla_x v dx dt \quad \text{and} \quad \ell(v) = \int_Q f v dx dt, \tag{2.3}$$

where the trial and test functions are defined by $H_{0,\underline{0}}^{1,0}(Q) = \{u \in L_2(Q) : \nabla_x u \in [L_2(Q)]^d, u = 0 \text{ on } \Sigma, \text{ and } u = 0 \text{ on } \Sigma_0\}$ and $H_{0,\overline{0}}^{1,1}(Q) = \{u \in L_2(Q) : \nabla_x u \in [L_2(Q)]^d, \partial_t u \in L_2(Q), u = 0 \text{ on } \Sigma, \text{ and } u = 0 \text{ on } \Sigma_T\}$. We denote the $\Sigma_T := \Omega(T)$, and $\nabla_x u = (\partial u/\partial x_1, \ldots, \partial u/\partial x_d)^\top$ denotes the gradient with respect to the spatial variables. The variational problem (2.2), including more general parabolic initial-boundary value problems, and other boundary conditions and more general elliptic parts, and nonlinear versions are known to have a unique weak solution, see e.g. [16, 31, 28].

## 2.2  B-Spline and NURBS Isogeometric Analysis

B-spline and NURBS based Galerkin methods has been studied and presented in the monograph [8]. Here, we present a briefly some information on the construction of the B-spline. Let us consider an interval $[0, 1]$, a vector $\Xi = \{0 = \xi_1, \ldots, \xi_{n+p+1} = 1\}$ with $p \geq 1$, and $n$ the number of basis functions. A *knot vector* a non-decreasing sequence of real numbers in the unit interval also called parameter domain $[0, 1]$. The knot vector $\Xi$ defines the uni-variate B-spline basis functions via the Cox -de Boor recursion formula

$$\widehat{B}_{i,0}(\xi) = \begin{cases} 1 & \text{if } \xi_i \leq \xi < \xi_{i+1}, \\ 0 & \text{else}, \end{cases}$$
$$\widehat{B}_{i,p}(\xi) = \frac{\xi - \xi_i}{\xi_{i+p} - \xi_i} \widehat{B}_{i,p-1}(\xi) + \frac{\xi_{i+p+1} - \xi}{\xi_{i+p+1} - \xi_{i+1}} \widehat{B}_{i+1,p-1}(\xi), \tag{2.4}$$

where a division by zero is defined to be zero. We note that a basis function of degree $p$ is $(p - m)$ times continuously differentiable across a knot value with the multiplicity $m$. If all internal knots have the multiplicity $m = 1$, then B-splines of degree $p$ are globally $C^{p-1}-$continuous.

In general for $(d+1)-$dimensional problems, the B-spline basis functions are tensor products of the univariate B-spline basis functions. Let $\Xi_\alpha = \{\xi_{1,\alpha}, \ldots, \xi_{n_\alpha+p_\alpha+1,\alpha}\}$ be the knot vectors for every direction $\alpha = 1, \ldots, d+1$. Let $\mathbf{i} := (i_1, \ldots, i_d)$, $\mathbf{p} := (p_1, \ldots, p_{d+1})$ and the set $\overline{\mathcal{I}} = \{\mathbf{i} = (i_1, \ldots, i_{d+1}) : i_\alpha = 1, 2, \ldots, n_\alpha; \ \alpha = 1, 2, \ldots, d+1\}$ be multi-indicies. Then the tensor product B-spline basis functions are defined by

$$\widehat{B}_{\mathbf{i},\mathbf{p}}(\xi) := \prod_{\alpha=1}^{d+1} \widehat{B}_{i_\alpha, p_\alpha}(\xi_\alpha), \tag{2.5}$$

where $\xi = (\xi_1, \ldots, \xi_{d+1}) \in \widehat{Q} = (0,1)^{d+1}$. The univariate and multivariate B-spline basis functions are defined in the parametric domain by means of the corresponding B-spline basis functions $\{\widehat{B}_{\mathbf{i},\mathbf{p}}\}_{\mathbf{i}\in\overline{\mathcal{I}}}$.

The distinct values $\xi_i, i = 1, \ldots, n$ of the knot vectors $\Xi$ provides a partition of $(0,1)^{d+1}$ creating a mesh $\widehat{\mathcal{K}}_h$ in the parameter domain where $\widehat{K}$ is a mesh element. The computational domain is described by means of a geometrical mapping $\mathbf{\Phi}$ such that $Q = \mathbf{\Phi}(\widehat{Q})$ and

$$\mathbf{\Phi}(\xi) := \sum_{\mathbf{i}\in\overline{\mathcal{I}}} C_{\mathbf{i}} \widehat{B}_{\mathbf{i},\mathbf{p}}(\xi), \tag{2.6}$$

where $C_{\mathbf{i}}$ are the control points. We define the basis functions in the computational domain by means of the geometrical mapping as $B_{\mathbf{i},\mathbf{p}} := \widehat{B}_{\mathbf{i},\mathbf{p}} \circ \mathbf{\Phi}^{-1}$ and the discrete function space by

$$\mathbb{V}_h = \mathrm{span}\{B_{\mathbf{i},\mathbf{p}} : \mathbf{i} \in \overline{\mathcal{I}}\}, \tag{2.7}$$

for the sake of the current work, we assume $\mathbb{V}_h \subset H^s(Q), s \geq 2$ i.e. continuously-differentiable $(C^1)$ basis functions on each patch.

In many practical applications, the computational domain $Q$ is decomposed into $N$ non-overlapping domain $Q_i$ called subdomains or patches denoted as $\mathcal{T}_h := \{Q_i\}_{i=1}^N$ such that $\overline{Q} = \bigcup_{i=1}^N \overline{Q}_i$ and $Q_i \cap Q_j = \emptyset$ for $i \neq j$. Each patch is the image of an associated geometrical mapping $\mathbf{\Phi}_i$ such that $\mathbf{\Phi}_i(\widehat{Q}) = Q_i, i = 1, \ldots, N$ see Fig. 1.

Let $F_{ij} = \partial Q_i \cap \partial Q_j, i \neq j$, denote the interior facets of two patches. The collection of all such interior facets is denoted by $\mathcal{F}_I$, and the collection of all Dirichlet facets $F_i = \partial Q_i \cap \partial Q$ is also denoted by $\mathcal{F}_D$. Furthermore, the collection of all internal and Dirichlet facets is denoted by $\mathcal{F} := \mathcal{F}_I \cup \mathcal{F}_D$.

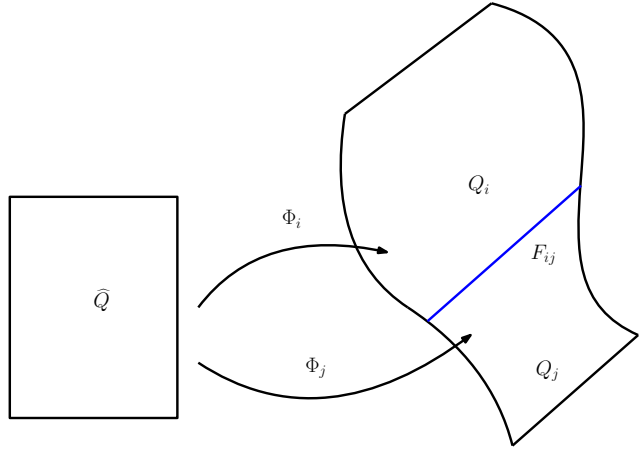

**Fig. 1.** A representation of the multi-patch space-time mapping

Also, we assume that for each patch $Q_i, i = 1, \ldots, N$, the underlying mesh $\mathcal{K}_{h,i}$ is quasi-uniform i.e.

$$h_K \leq h_i \leq C_q h_K, \quad \text{for all} \quad K \in \mathcal{K}_{h,i}, \quad i = 1, \ldots, N, \qquad (2.8)$$

where $C_q \geq 1$ and $h_i = \max\{h_K, K \in \mathcal{K}_{h,i}\}$ is the mesh size of $Q_i$ and $h_K$ is the diameter of of the mesh element $K$.

## 3    Stable Space-Time Discontinuous Galerkin Isogeometric Analysis

For a non-overlapping decomposition of the space-time computational domain $Q$ into $N$ patches, we assign non-negative integers $s_i$ to each patch $Q_i$, and collect them in the vector $\mathbf{s} = \{s_1, \ldots, s_N\}$. Let us now define the broken Sobolev space

$$H^{\mathbf{s}}(Q, \mathcal{T}_h) := \{v \in L_2(Q) : v|_{Q_i} \in H^{s_i}(Q_i), \ \forall \, i = 1, \ldots, N\}, \qquad (3.1)$$

and equip it with a broken Sobolev norm and semi-norm

$$\|v\|_{H^{\mathbf{s}}(Q, \mathcal{T}_h)} := \left( \sum_{i=1}^{N} \|v\|^2_{H^{s_i}(Q_i)} \right)^{1/2} \quad \text{and} \quad |v|_{H^{\mathbf{s}}(Q, \mathcal{T}_h)} := \left( \sum_{i=1}^{N} |v|^2_{H^{s_i}(Q_i)} \right)^{1/2},$$
$$(3.2)$$

respectively. Next, we define the discrete space

$$V_{0h} := \{v_h \in L_2(Q) : v_h|_{Q_i} \in \mathbb{V}_{h,i}, i = 1, \ldots, N \quad \text{and} \quad v_h = 0 \quad \Sigma_0 \cap \Sigma\}. \tag{3.3}$$

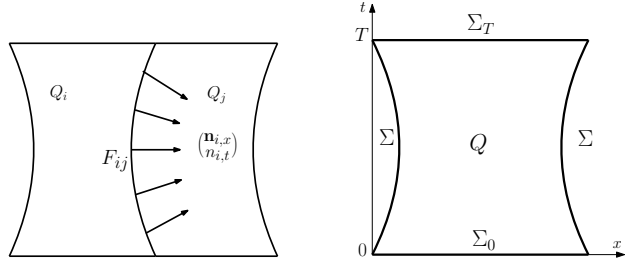

**Fig. 2.** A representation of the interface $F_{ij}$ of two neighboring patches $Q_i$ and $Q_j$ and its associated normal vector $\mathbf{n}_i = (\mathbf{n}_{i,x}, n_{i,t})$ on the left. The underlying physical mesh (blue) for the multi-patch on the right.

**Definition 1.** *Let $F_{ij} \in \mathcal{F}_I$ be an interior facet with the outer unit normal vector $\mathbf{n}_i = (\mathbf{n}_{i,x}, n_{i,t})^\top \in \mathbb{R}^{d+1}$ with respect to $Q_i \subset \mathbb{R}^{d+1}$ with $d \in \{1, 2, 3\}$. For a given, sufficiently smooth scalar function $v$, we will denote the restriction of the function $v$ to $Q_i$ and $Q_j$ by $v_i$ and $v_j$, respectively. We define the jumps across the facets by*

$$\llbracket v \rrbracket := \begin{cases} v_i \mathbf{n}_i + v_j \mathbf{n}_j & on \quad F_{ij} \in \mathcal{F}_I, \\ v_i \mathbf{n}_i & on \quad F_i \in \mathcal{F} \setminus \mathcal{F}_I. \end{cases}$$

*The jump in space direction is given by*

$$\llbracket v \rrbracket_x := \begin{cases} v_i \mathbf{n}_{i,x} + v_j \mathbf{n}_{j,x} & on \quad F_{ij} \in \mathcal{F}_I, \\ v_i \mathbf{n}_{i,x} & on \quad F_i \in \mathcal{F} \setminus \mathcal{F}_I, \end{cases}$$

*whereas the jump in time direction is defined by*

$$\llbracket v \rrbracket_t := \begin{cases} v_i n_{i,t} + v_j n_{j,t} & on \quad F_{ij} \in \mathcal{F}_I, \\ v_i n_{i,t} & on \quad F_i \in \mathcal{F} \setminus \mathcal{F}_I. \end{cases}$$

*The average of function $v$ is nothing but*

$$\{v\} := \begin{cases} \dfrac{1}{2}(v_i + v_j) & on \quad F_{ij} \in \mathcal{F}_I, \\ v_i & on \quad F_i \in \mathcal{F} \setminus \mathcal{F}_I. \end{cases}$$

*For the interior facets $F_{ij} \in \mathcal{F}_I$, the upwind in time direction is given by*

$$\{v\}^{up} := \begin{cases} v_i & for \quad n_{i,t} \geq 0, \\ v_j & for \quad n_{i,t} < 0, \end{cases} \tag{3.4}$$

*whereas the downwind in time direction is given by*

$$\{v\}^{down} := \begin{cases} v_j & for \quad n_{i,t} \geq 0, \\ v_i & for \quad n_{i,t} < 0. \end{cases} \tag{3.5}$$

The derivation of the variational scheme requires some formulas stated in the next lemma.

**Lemma 1.** *Let $F_{ij} \in \mathcal{F}_I$ be an interface, and let $u$ and $v$ be sufficiently smooth functions. Then the following formulas hold:*

$$[\![uv]\!]_x = \{u\}[\![v]\!]_x + \{v\}[\![u]\!]_x, \tag{3.6}$$

$$[\![uv]\!]_t = \{u\}^{up}[\![v]\!]_t + \{v\}^{down}[\![u]\!]_t. \tag{3.7}$$

*Proof.* See e.g. [21, Lemma 5.17]. $\qquad\square$

We will begin with the derivation of the space-time multipatch dG-IgA scheme as follows: We multiply the model problem (2.1) by a time-upwind test function $v_h + \theta h \partial_t v_h$ with an arbitrary $v_h \in V_{0h}$ with a positive constant $\theta$ and sum over each patch

$$\sum_{i=1}^{N} \int_{Q_i} f(v_h + \theta_i h_i \partial_t v_h) \, dx dt$$

$$= \sum_{i=1}^{N} \int_{Q_i} \left( \partial_t u(v_h + \theta_i h_i \partial_t v_h) - \Delta u(v_h + \theta_i h_i \partial_t v_h) \right) dx dt. \tag{3.8}$$

Concerning the Laplacian term $\Delta u$, we sum over each patch and apply an integration by parts with respect to the spatial direction leading to

$$-\sum_{i=1}^{N} \int_{Q_i} \Delta u(v_h + \theta_i h_i \partial_t v_h) \, dx dt$$

$$= \sum_{i=1}^{N} \int_{Q_i} \left( \nabla_x u \cdot \nabla_x v_h + \theta_i h_i \nabla_x u \cdot \nabla_x \partial_t v_h \right) dx dt$$

$$- \sum_{i=1}^{N} \int_{\partial Q_i} \mathbf{n}_{i,x} \cdot \nabla_x u(v_h + \theta_i h_i \partial_t v_h) \, ds. \tag{3.9}$$

Integrating by parts the first term with respect to the time yields

$$= \sum_{i=1}^{N} \int_{Q_i} \left( \nabla_x u \cdot \nabla_x v_h - \theta_i h_i \partial_t \nabla_x u \cdot \nabla_x v_h \right) dx dt$$

$$+ \theta_i h_i \sum_{i=1}^{N} \int_{\partial Q_i} n_{i,t} \nabla_x u \cdot \nabla_x v_h \, ds$$

$$- \sum_{i=1}^{N} \int_{\partial Q_i} \left( \mathbf{n}_{i,x} \cdot \nabla_x u (v_h + \theta_i h_i \partial_t v_h) \right) ds. \qquad (3.10)$$

By considering the boundary terms in (3.10), we rewrite the sum over the boundary integrals as

$$\theta_i h_i \sum_{i=1}^{N} \int_{\partial Q_i} \left( n_{i,t} \nabla_x u \cdot \nabla_x v_h \right) ds = \theta_i h_i \sum_{F_{ij} \in \mathcal{F}_I} \int_{F_{ij}} [\![ \nabla_x u \cdot \nabla_x v_h ]\!]_t \, ds$$

$$+ \theta_i h_i \sum_{F_i \in \mathcal{F}_0 \cup \mathcal{F}_T \cup \mathcal{F}_D} \int_{F_i} n_{i,t} \nabla_x u \cdot \nabla_x v_h \, ds, \qquad (3.11)$$

and

$$\sum_{i=1}^{N} \int_{\partial Q_i} \left( \mathbf{n}_{i,x} \cdot \nabla_x u (v_h + \theta_i h_i \partial_t v_h) \right) ds$$

$$= \sum_{F_{ij} \in \mathcal{F}_I} \int_{F_{ij}} [\![ \nabla_x u v_h ]\!]_x \, ds + \theta_i h_i \sum_{F_{ij} \in \mathcal{F}_I} \int_{F_{ij}} [\![ \nabla_x u \partial_t v_h ]\!]_x \, ds$$

$$+ \sum_{F_i \in \mathcal{F}_0 \cup \mathcal{F}_T \cup \mathcal{F}_D} \int_{F_i} \mathbf{n}_{i,x} \cdot \nabla_x u (v_h + \theta_i h_i \partial_t v_h) \, ds. \qquad (3.12)$$

Since $v_h = 0$ on the facets $F_i \in \mathcal{F}_0$, the term (3.11) yields

$$\theta_i h_i \sum_{i=1}^{N} \int_{\partial Q_i} \left( n_{i,t} \nabla_x u \cdot \nabla_x v_h \right) ds = \theta_i h_i \sum_{F_{ij} \in \mathcal{F}_I} \int_{F_{ij}} [\![ \nabla_x u \cdot \nabla_x v_h ]\!]_t \, ds$$

$$+ \theta_i h_i \sum_{F_i \in \mathcal{F}_T \cup \mathcal{F}_D} \int_{F_i} n_{i,t} \nabla_x u \cdot \nabla_x v_h \, ds. \qquad (3.13)$$

The normal vector in the space direction $\mathbf{n}_{i,x}$ is zero for $F_i \in \mathcal{F}_0$ and $F_i \in \mathcal{F}_T$, thus (3.12) can be rewritten as

$$\sum_{i=1}^{N} \int_{\partial Q_i} \left( \mathbf{n}_{i,x} \cdot \nabla_x u(v_h + \theta_i h_i \partial_t v_h) \right) ds = \sum_{F_{ij} \in \mathcal{F}_I} \int_{F_{ij}} [\![ \nabla_x u v_h ]\!]_x \, ds$$

$$+ \theta_i h_i \sum_{F_{ij} \in \mathcal{F}_I} \int_{F_{ij}} [\![ \nabla_x u \partial_t v_h ]\!]_x \, ds + \theta_i h_i \sum_{F_i \in \mathcal{F}_D} \int_{F_i} \mathbf{n}_{i,x} \cdot \nabla_x u \partial_t v_h \, ds.$$

$$(3.14)$$

Substituting (3.13) and (3.14) into (3.10) and rearranging the terms yields

$$-\sum_{i=1}^{N} \int_{Q_i} \Delta u (v_h + \theta_i h_i \partial_t v_h) \, dxdt = \sum_{i=1}^{N} \int_{Q_i} \left( \nabla_x u \cdot \nabla_x v_h - \theta_i h_i \partial_t \nabla_x u \cdot \nabla_x v_h \right) dxdt$$

$$+ \theta_i h_i \sum_{F_i \in \mathcal{F}_T} \int_{F_i} n_{i,t} \nabla_x u \cdot \nabla_x v_h \, ds + \theta_i h_i \sum_{F_{ij} \in \mathcal{F}_I} \int_{F_{ij}} [\![ \nabla_x u \cdot \nabla_x v_h ]\!]_t \, ds$$

$$- \sum_{F_{ij} \in \mathcal{F}_I} \int_{F_{ij}} \left( [\![ \nabla_x u v_h ]\!]_x + \theta_i h_i [\![ \nabla_x u \partial_t v_h ]\!]_x \right) ds$$

$$+ \theta_i h_i \sum_{F_i \in \mathcal{F}_D} \int_{F_i} \left( n_{i,t} \nabla_x u \cdot \nabla_x v_h - \mathbf{n}_{i,x} \cdot \nabla_x u \partial_t v_h \right) ds.$$

Now, by using Lemma 1, which also holds for vector functions, we can rewrite the terms on the interior facets $F_{ij}$ to obtain

$$= \sum_{i=1}^{N} \int_{Q_i} \left( \nabla_x u \cdot \nabla_x v_h - \theta_i h_i \partial_t \nabla_x u \cdot \nabla_x v_h \right) dxdt + \theta_i h_i \sum_{F_i \in \mathcal{F}_T} \int_{F_i} \nabla_x u \cdot \nabla_x v_h \, ds$$

$$+ \theta_i h_i \sum_{F_{ij} \in \mathcal{F}_I} \int_{F_{ij}} \left( \{ \nabla_x u \}^{up} [\![ \nabla_x v_h ]\!]_t + \{ \nabla_x v_h \}^{down} [\![ \nabla_x u ]\!]_t \right) ds$$

$$- \sum_{F_{ij} \in \mathcal{F}_I} \int_{F_{ij}} \left( \{ \nabla_x u \} [\![ v_h ]\!]_x \, ds + [\![ \nabla_x u ]\!]_x \{ v_h \} \right) ds$$

$$- \theta_i h_i \sum_{F_{ij} \in \mathcal{F}_I} \int_{F_{ij}} \left( \{ \nabla_x u \} [\![ \partial_t v_h ]\!]_x \, ds + [\![ \nabla_x u ]\!]_x \{ \partial_t v_h \} \right) ds$$

$$+ \theta_i h_i \sum_{F_i \in \mathcal{F}_D} \int_{F_i} \left( n_{i,t} \nabla_x u \cdot \nabla_x v_h - \mathbf{n}_{i,x} \cdot \nabla_x u \partial_t v_h \right) ds.$$

The jumps of the derivatives of the solution $u$ are zero, i.e., $[\![\nabla_x u]\!]_x = 0$ and $[\![\nabla_x u]\!]_t = 0$ yielding

$$
\begin{aligned}
&= \sum_{i=1}^{N} \int_{Q_i} \left( \nabla_x u \cdot \nabla_x v_h - \theta_i h_i \partial_t \nabla_x u \cdot \nabla_x v_h \right) dx dt + \theta_i h_i \sum_{F_i \in \mathcal{F}_T} \int_{F_i} \nabla_x u \cdot \nabla_x v_h \, ds \\
&+ \theta_i h_i \sum_{F_{ij} \in \mathcal{F}_I} \int_{F_{ij}} \{\nabla_x u\}^{up} [\![\nabla_x v_h]\!]_t \, ds - \sum_{F_{ij} \in \mathcal{F}_I} \int_{F_{ij}} \{\nabla_x u\} [\![v_h]\!]_x \, ds \\
&- \theta_i h_i \sum_{F_{ij} \in \mathcal{F}_I} \int_{F_{ij}} \{\nabla_x u\} [\![\partial_t v_h]\!]_x \, ds + \theta_i h_i \sum_{F_i \in \mathcal{F}_D} \int_{F_i} \left( n_{i,t} \nabla_x u \cdot \nabla_x v_h - \mathbf{n}_{i,x} \cdot \nabla_x u \partial_t v_h \right) ds.
\end{aligned}
$$

Also, for the exact solution $u$, we have $[\![u]\!]_x = 0$ and $[\![\partial_t u]\!]_x = 0$ on the interior facets $F_{ij} \in \mathcal{F}_I$. Therefore, we add the following consistent terms

$$
\sum_{F_{ij} \in \mathcal{F}_I} \int_{F_{ij}} \{\nabla_x v_h\} [\![u]\!]_x \, ds + \theta_i h_i \sum_{F_{ij} \in \mathcal{F}_I} \int_{F_{ij}} \{\nabla_x v_h\} [\![\partial_t u]\!]_x \, ds.
$$

Furthermore, we penalize the interior facets by adding consistency terms that are zero to the bilinear form

$$
\sum_{F_{ij} \in \mathcal{F}_I} \int_{F_{ij}} \frac{\delta_1}{h_i} [\![u]\!]_x [\![v_h]\!]_x \, ds + \sum_{F_{ij} \in \mathcal{F}_I} \int_{F_{ij}} \delta_2 \theta_i h_i [\![\partial_t u]\!]_x [\![\partial_t v_h]\!]_x \, ds,
$$

where the penalty parameters $\delta_1$ and $\delta_2$ are positive constants that will be determined later. Considering the terms on the Dirichlet facets

$$
\begin{aligned}
\sum_{F_i \in \mathcal{F}_D} \int_{F_i} & \left( n_{i,t} \nabla_x u \cdot \nabla_x v_h - \mathbf{n}_{i,x} \cdot \nabla_x u \partial_t v_h \right) ds \\
&= \sum_{F_i \in \mathcal{F}_D} \int_{F_i} \nabla_x u \cdot \left( n_{i,t} \nabla_x v_h - \mathbf{n}_{i,x} \partial_t v_h \right) ds \qquad (3.15)
\end{aligned}
$$

Since $(n_{i,t}\nabla_x v_h - \mathbf{n}_{i,x}\partial_t v_h)$ is the tangential derivative of $v_h$ and $v_h = 0$ on the facets $F_i \in \mathcal{F}_D$. Thus, we can proceed as follows

$$= \sum_{i=1}^{N} \int_{Q_i} \left( \nabla_x u \cdot \nabla_x v_h - \theta_i h_i \partial_t \nabla_x u \cdot \nabla_x v_h \right) dxdt + \theta_i h_i \sum_{F_i \in \mathcal{F}_T} \int_{F_i} \nabla_x u \cdot \nabla_x v_h \, ds$$

$$+ \theta_i h_i \sum_{F_{ij} \in \mathcal{F}_I} \int_{F_{ij}} \{\nabla_x u\}^{up} [\![\nabla_x v_h]\!]_t \, ds + \sum_{F_{ij} \in \mathcal{F}_I} \int_{F_{ij}} \left( - \{\nabla_x u\} [\![v_h]\!]_x + \{\nabla_x v_h\} [\![u]\!]_x \right) ds$$

$$+ \theta_i h_i \sum_{F_{ij} \in \mathcal{F}_I} \int_{F_{ij}} \left( - \{\nabla_x u\} [\![\partial_t v_h]\!]_x \, ds + \{\nabla_x v_h\} [\![\partial_t u]\!]_x \right) ds$$

$$+ \sum_{F_{ij} \in \mathcal{F}_I} \int_{F_{ij}} \left( \frac{\delta_1}{h_i} [\![u]\!]_x [\![v_h]\!]_x + \delta_2 \theta_i h_i [\![\partial_t u]\!]_x [\![\partial_t v_h]\!]_x \right) ds.$$

Next, considering the time derivative term $\partial_t u$, we apply an integration by part in time as follows

$$\sum_{i=1}^{N} \int_{Q_i} \partial_t u v_h \, dxdt = -\sum_{i=1}^{N} \int_{Q_i} u \partial_t v_h \, dxdt + \sum_{i=1}^{N} \int_{\partial Q_i} n_{i,t} u v_h \, ds. \quad (3.16)$$

Rewriting the sum over the boundary integral terms in (3.16), we obtain the following identity

$$\sum_{i=1}^{N} \int_{\partial Q_i} n_{i,t} u v_h \, ds = \sum_{F_{ij} \in \mathcal{F}_I} \int_{F_{ij}} [\![u v_h]\!]_t \, ds + \sum_{F_i \in \mathcal{F}_0 \cup \mathcal{F}_T \cup \mathcal{F}_D} \int_{F_i} n_{i,t} u v_h \, ds. \quad (3.17)$$

Using (3.7) of Lemma 1 and the fact that $v_h = 0$ on the facets $F_i \in \mathcal{F}_D \cup \mathcal{F}_0$ and $n_{i,t} = 1$ on $F_i \in \mathcal{F}_T$, the boundary integral term yields

$$\sum_{i=1}^{N} \int_{\partial Q_i} n_{i,t} u v_h \, ds = \sum_{F_{ij} \in \mathcal{F}_I} \int_{F_{ij}} \left( \{u\}^{up} [\![v_h]\!]_t + \{v_h\}^{down} [\![u]\!]_t \right) ds$$

$$+ \sum_{F_i \in \mathcal{F}_T} \int_{F_i} u v_h \, ds. \quad (3.18)$$

For a smooth exact solution $u$, the jump in time direction is zero, i.e. $[\![u]\!]_t = 0$. Thus, we obtain

$$\sum_{i=1}^{N} \int_{\partial Q_i} n_{i,t} u v_h \, ds = \sum_{F_{ij} \in \mathcal{F}_I} \int_{F_{ij}} \{u\}^{up} [\![v_h]\!]_t \, ds + \sum_{F_i \in \mathcal{F}_T} \int_{F_i} u v_h \, ds. \quad (3.19)$$

In summary, we have shown that a weak solution $u \in H^2(Q, \mathcal{T}_h)$ of our variational problem (2.2) fulfills the space-time dG-IgA variational identity

$$a_h(u, v_h) = \ell_h(v_h) \quad \forall v_h \in V_{0h}, \tag{3.20}$$

where the bilinear form is represented by

$$a_h(u, v_h) = \sum_{i=1}^{N} \int_{Q_i} \left( -u \partial_t v_h + \theta_i h_i \partial_t u \partial_t v_h + \nabla_x u \cdot \nabla_x v_h - \theta_i h_i \partial_t \nabla_x u \cdot \nabla_x v_h \right) dxdt$$

$$+ \sum_{F_{ij} \in \mathcal{F}_I} \int_{F_{ij}} \{u\}^{up} [\![v_h]\!]_t \, ds + \sum_{F_i \in \mathcal{F}_T} \int_{F_i} uv_h \, ds + \theta_i h_i \sum_{F_i \in \mathcal{F}_T} \int_{F_i} \nabla_x u \cdot \nabla_x v_h \, ds$$

$$+ \theta_i h_i \sum_{F_{ij} \in \mathcal{F}_I} \int_{F_{ij}} \{\nabla_x u\}^{up} [\![\nabla_x v_h]\!]_t \, ds + \sum_{F_{ij} \in \mathcal{F}_I} \int_{F_{ij}} \left( -\{\nabla_x u\} [\![v_h]\!]_x + \{\nabla_x v_h\} [\![u]\!]_x \right) ds$$

$$+ \theta_i h_i \sum_{F_{ij} \in \mathcal{F}_I} \int_{F_{ij}} \left( -\{\nabla_x u\} [\![\partial_t v_h]\!]_x \, ds + \{\nabla_x v_h\} [\![\partial_t u]\!]_x \right) ds$$

$$+ \sum_{F_{ij} \in \mathcal{F}_I} \int_{F_{ij}} \left( \frac{\delta_1}{h_i} [\![u]\!]_x [\![v_h]\!]_x + \delta_2 \theta_i h_i [\![\partial_t u]\!]_x [\![\partial_t v_h]\!]_x \right) ds, \tag{3.21}$$

and the linear form

$$\ell(v) = \sum_{i=1}^{N} \int_{Q_i} f \left( v_h + \theta_i h_i \partial_t v_h \right) dxdt. \tag{3.22}$$

*Remark 1.* We note that the penalty parameters $\delta_1$ and $\delta_2$ are known to depend on the NURBS degree $p$ and spatial dimension of the computational domain $d$. The space-time dG-IgA is a non-symmetric discretization scheme.

Next, we show that the bilinear form (3.21) is coercive with respect to the discrete norm

$$\|v_h\|_h^2 = \sum_{i=1}^{N} \|\nabla_x v_h\|_{L_2(Q_i)}^2 + \theta_i h_i \|\partial_t v_h\|_{L_2(Q_i)}^2 + \frac{1}{2} \sum_{F_i \in \mathcal{F}_T} \|v_h\|_{L_2(F_i)}^2$$

$$+ \frac{\theta_i h_i}{2} \sum_{F_i \in \mathcal{F}_T} \|\nabla_x v_h\|_{L_2(F_i)}^2 + \frac{1}{2} \sum_{F_{ij} \in \mathcal{F}_I} \|[\![v_h]\!]_t\|_{L_2(F_{ij})}^2 + \frac{\theta_i h_i}{2} \sum_{F_{ij} \in \mathcal{F}_I} \|[\![\nabla_x v_h]\!]_t\|_{L_2(F_{ij})}^2$$

$$+ \sum_{F_{ij} \in \mathcal{F}_I} \left( \frac{\delta_1}{h_i} \|[\![v_h]\!]_x\|_{L_2(F_{ij})}^2 + \delta_2 \theta_i h_i \|[\![\partial_t v_h]\!]_x\|_{L_2(F_{ij})}^2 \right). \tag{3.23}$$

*Remark 2.* The space-time dG norm $\|\cdot\|_h$ from (3.23) is a norm on $V_{0h}$. Indeed, if $\|v_h\|_h = 0$ for some function $v_h \in V_{0h}$, then $\nabla_x v_h = 0$ and $\partial_t v_h = 0$ in each subdomain $Q_i$. This means that the function $v_h$ is a constant on each patch $Q_i$, $i = 1, \ldots, N$. Furthermore, $\|v_h\|_h = 0$ yields that the jumps $[\![v_h]\!]_x, [\![v_h]\!]_t, [\![\partial_t v_h]\!]_x$ and $[\![\nabla_x v_h]\!]_t$ of $v_h$ across the internal facets $F_{ij} \in \mathcal{F}_I$ are zero, i.e., $v_h$ is constant in $\overline{Q}$. Finally, $v_h \in V_{0h}$ implies that $v_h$ is zero on the facets $F_i \in \mathcal{F}_D$ and $F_i \in \mathcal{F}_0$, i.e. this constant must be zero. Therefore, $v_h = 0$ in the whole space-time computational domain $Q$. The other norm axioms (homogeneity and triangle inequality) are easily verified.

The following lemmata are required to prove the coercivity of the bilinear form.

**Lemma 2.** *Let $F_{ij} \in \mathcal{F}_I$ be an interior facet and $v \in V_{0h}$ be a function. Then the following identity holds:*

$$\{v\}^{up}[\![v]\!]_t - \frac{1}{2}[\![v^2]\!]_t = \frac{1}{2}|n_{i,t}|[\![v]\!]^2. \tag{3.24}$$

*Proof.* See e.g. [21, Lemma 5.19].

**Lemma 3.** *Let $K$ be an arbitrary mesh element from $\mathcal{K}_{h,i}$; the underlying mesh of $Q_i, i = 1, \ldots, N$. Then the inverse inequalities*

$$\|\nabla v\|^2_{L_2(Q_i)} \leq C_{inv,1,u} h_i^{-2} \|v\|^2_{L_2(Q_i)}, \tag{3.25}$$

*and*

$$\|v\|^2_{L_2(\partial Q_i)} \leq C_{inv,0,u} h_i^{-1} \|v\|^2_{L_2(Q_i)}, \tag{3.26}$$

*hold for all $v \in V_h$, where $C_{inv,0,u}$ and $C_{inv,1,u}$ are positive constants, which are independent of $h_i$ and $Q_i$.*

*Proof.* See e.g. [21, Chapter 2].

**Lemma 4.** *Let $\theta_i > 0$ be sufficiently small such that $\theta_i \leq (C_{inv,0,u} C_q)^{-1}$ where $C_{inv,0,u}$ is chosen as in Lemma 3 and $C_q$ satisfies (2.8). Let $\delta_1$ and $\delta_2$ be some positive constants, then the discrete bilinear form $a_h(\cdot, \cdot) : V_{0h} \times V_{0h} \to \mathbb{R}$, defined by (3.21), is $V_{0h}-$coercive with respect to the norm $\|\cdot\|_h$, i.e.*

$$a_h(v_h, v_h) \geq \mu_c \|v_h\|^2_h, \quad \forall v_h \in V_{0h}, \tag{3.27}$$

*where $\mu_c = 1/2$.*

*Proof.* By using the definition of the bilinear form (3.21), we proceed as follows

$$a_h(v_h, v_h) = \sum_{i=1}^{N} -\frac{1}{2} \int_{Q_i} \partial_t v_h^2 \, dx dt + \theta_i h_i \|\partial_t v_h\|_{L_2(Q_i)}^2 + \|\nabla_x v_h\|_{L_2(Q_i)}^2$$

$$-\frac{\theta_i h_i}{2} \int_{Q_i} \partial_t |\nabla_x v_h|^2 \, dx dt + \sum_{F_{ij} \in \mathcal{F}_I} \int_{F_{ij}} \{v_h\}^{up} [\![v_h]\!]_t \, ds + \sum_{F_i \in \mathcal{F}_T} \int_{F_i} v_h^2 \, ds$$

$$+ \theta_i h_i \sum_{F_i \in \mathcal{F}_T} \|\nabla_x v_h\|_{L_2(F_i)}^2 + \theta_i h_i \sum_{F_{ij} \in \mathcal{F}_I} \int_{F_{ij}} \{\nabla_x v_h\}^{up} [\![\nabla_x v_h]\!]_t \, ds$$

$$+ \sum_{F_{ij} \in \mathcal{F}_I} \left( \frac{\delta_1}{h} \| [\![v_h]\!]_x \|_{L_2(F_{ij})}^2 + \delta_2 \theta_i h_i \| [\![\partial_t v_h]\!]_x \|_{L_2(F_{ij})}^2 \right)$$

$$= \sum_{i=1}^{N} \left( -\frac{1}{2} \int_{\partial Q_i} n_{i,t} v_h^2 \, ds + \theta_i h_i \|\partial_t v_h\|_{L_2(Q_i)}^2 + \|\nabla_x v_h\|_{L_2(Q_i)}^2 \right.$$

$$\left. -\frac{\theta_i h_i}{2} \int_{\partial Q_i} n_{i,t} |\nabla_x v_h|^2 \, ds \right) + \sum_{F_{ij} \in \mathcal{F}_I} \int_{F_{ij}} \{v_h\}^{up} [\![v_h]\!]_t \, ds + \sum_{F_i \in \mathcal{F}_T} \int_{F_i} v_h^2 \, ds$$

$$+ \theta_i h_i \sum_{F_i \in \mathcal{F}_T} \|\nabla_x v_h\|_{L_2(F_i)}^2 + \theta_i h_i \sum_{F_{ij} \in \mathcal{F}_I} \int_{F_{ij}} \{\nabla_x v_h\}^{up} [\![\nabla_x v_h]\!]_t \, ds$$

$$+ \sum_{F_{ij} \in \mathcal{F}_I} \left( \frac{\delta_1}{h} \| [\![v_h]\!]_x \|_{L_2(F_{ij})}^2 + \delta_2 \theta_i h_i \| [\![\partial_t v_h]\!]_x \|_{L_2(F_{ij})}^2 \right).$$

Using Gauss' theorem with $v_h = 0$ on the facet $F_i \in \mathcal{F}_0$, we proceed as follows

$$= \sum_{i=1}^{N} \left( \theta_i h_i \|\partial_t v_h\|_{L_2(Q_i)}^2 + \|\nabla_x v_h\|_{L_2(Q_i)}^2 \right) - \frac{1}{2} \sum_{F_i \in \mathcal{F}_T} \int_{F_i} v_h^2 \, ds$$

$$-\frac{1}{2} \sum_{F_{ij} \in \mathcal{F}_I} \int_{F_{ij}} [\![v_h^2]\!]_t \, ds + \sum_{F_{ij} \in \mathcal{F}_I} \int_{F_{ij}} \{v_h\}^{up} [\![v_h]\!]_t \, ds + \sum_{F_i \in \mathcal{F}_T} \int_{F_i} v_h^2 \, ds$$

$$-\frac{\theta_i h_i}{2} \sum_{F_{ij} \in \mathcal{F}_I} \int_{F_{ij}} [\![|\nabla_x v_h|^2]\!]_t \, ds - \frac{\theta_i h_i}{2} \sum_{F_i \in \mathcal{F}_T} \|\nabla_x v_h\|_{L_2(F_i)}^2 - \frac{\theta_i h_i}{2} \sum_{F_i \in \mathcal{F}_D} \int_{F_i} n_{i,t} |\nabla_x v_h|^2 \, ds$$

$$+ \theta_i h_i \sum_{F_i \in \mathcal{F}_T} \|\nabla_x v_h\|_{L_2(F_i)}^2 + \theta_i h_i \sum_{F_{ij} \in \mathcal{F}_I} \int_{F_{ij}} \{\nabla_x v_h\}^{up} [\![\nabla_x v_h]\!]_t \, ds$$

$$+ \sum_{F_{ij} \in \mathcal{F}_I} \left( \frac{\delta_1}{h_i} \| [\![v_h]\!]_x \|_{L_2(F_{ij})}^2 + \delta_2 \theta_i h_i \| [\![\partial_t v_h]\!]_x \|_{L_2(F_{ij})}^2 \right)$$

$$= \sum_{i=1}^{N} \left( \theta_i h_i \|\partial_t v_h\|_{L_2(Q_i)}^2 + \|\nabla_x v_h\|_{L_2(Q_i)}^2 \right) + \frac{1}{2} \sum_{F_i \in \mathcal{F}_T} \|v_h\|_{L_2(F_i)}^2 + \frac{\theta_i h_i}{2} \sum_{F_i \in \mathcal{F}_T} \|\nabla_x v_h\|_{L_2(F_i)}^2$$

$$+ \sum_{F_{ij} \in \mathcal{F}_I} \int_{F_{ij}} \left( \{v_h\}^{up} [\![v_h]\!]_t - \frac{1}{2} [\![v_h^2]\!]_t \right) ds - \frac{\theta_i h_i}{2} \sum_{F_i \in \mathcal{F}_D} \int_{F_i} n_{i,t} |\nabla_x v_h|^2 \, ds$$

$$+ \theta_i h_i \sum_{F_{ij} \in \mathcal{F}_I} \int_{F_{ij}} \left( \{\nabla_x v_h\}^{up} [\![\nabla_x v_h]\!]_t - \frac{1}{2} [\![|\nabla_x v_h|^2]\!]_t \right) ds$$

$$+ \sum_{F_{ij} \in \mathcal{F}_I} \frac{\delta_1}{h_i} \|[\![v_h]\!]_x\|_{L_2(F_{ij})}^2 + \sum_{F_{ij} \in \mathcal{F}_I} \delta_2 \theta_i h_i \|[\![\partial_t v_h]\!]_x\|_{L_2(F_{ij})}^2.$$

Now, by using Lemma 2 together with $|n_{i,t}| \geq |n_{i,t}|^2$, we sum over $Q_i$ as follows

$$= \sum_{i=1}^{N} \left( \theta_i h_i \|\partial_t v_h\|_{L_2(Q_i)}^2 + \|\nabla_x v_h\|_{L_2(Q_i)}^2 \right) + \frac{1}{2} \sum_{F_i \in \mathcal{F}_T} \|v_h\|_{L_2(F_i)}^2$$

$$+ \frac{\theta_i h_i}{2} \sum_{F_i \in \mathcal{F}_T} \|\nabla_x v_h\|_{L_2(F_i)}^2 + \frac{1}{2} \sum_{F_{ij} \in \mathcal{F}_I} \int_{F_{ij}} |n_{i,t}| \left( [\![v_h]\!] \right)^2 ds$$

$$+ \frac{\theta_i h_i}{2} \sum_{F_{ij} \in \mathcal{F}_I} \int_{F_{ij}} |n_{i,t}| \left( [\![|\nabla_x v_h|]\!] \right)^2 ds - \frac{\theta_i h_i}{2} \sum_{F_i \in \mathcal{F}_D} \int_{F_i} n_{i,t} |\nabla_x v_h|^2 \, ds$$

$$+ \sum_{F_{ij} \in \mathcal{F}_I} \left( \frac{\delta_1}{h} \|[\![v_h]\!]_x\|_{L_2(F_{ij})}^2 + \delta_2 \theta_i h_i \|[\![\partial_t v_h]\!]_x\|_{L_2(F_{ij})}^2 \right)$$

$$\geq \sum_{i=1}^{N} \left( \theta_i h_i \|\partial_t v_h\|_{L_2(Q_i)}^2 + \|\nabla_x v_h\|_{L_2(Q_i)}^2 \right) + \frac{1}{2} \sum_{F_i \in \mathcal{F}_T} \|v_h\|_{L_2(F_i)}^2 + \frac{\theta_i h_i}{2} \sum_{F_i \in \mathcal{F}_T} \|\nabla_x v_h\|_{L_2(F_i)}^2$$

$$+ \frac{1}{2} \sum_{F_{ij} \in \mathcal{F}_I} \|[\![v_h]\!]_t\|_{L_2(F_{ij})}^2 + \frac{\theta_i h_i}{2} \sum_{F_{ij} \in \mathcal{F}_I} \|[\![\nabla_x v_h]\!]_t\|_{L_2(F_{ij})}^2 - \frac{\theta_i h_i}{2} \sum_{F_i \in \mathcal{F}_D} \|\nabla_x v_h\|_{L_2(F_i)}^2$$

$$+ \sum_{F_{ij} \in \mathcal{F}_I} \left( \frac{\delta_1}{h} \|[\![v_h]\!]_x\|_{L_2(F_{ij})}^2 + \sum_{F_{ij} \in \mathcal{F}_I} \delta_2 \theta_i h_i \|[\![\partial_t v_h]\!]_x\|_{L_2(F_{ij})}^2 \right)$$

$$= \|v_h\|_h^2 - \frac{\theta_i h_i}{2} \sum_{F_i \in \mathcal{F}_D} \|\nabla_x v_h\|_{L_2(F_i)}^2.$$

Next, we estimate the term on the facet $F_i \in \mathcal{F}_D$ by using the quasi-uniformity property (2.8) and the inverse inequality (3.26) to arrive at

$$(\theta_i h_i)/2 \sum_{F_i \in \mathcal{F}_D} \|\nabla_x v_h\|_{L_2(F_i)}^2 \leq (\theta_i C_q C_{inv,0,u})/2 \sum_{i=1}^{N} \|\nabla_x v_h\|_{L_2(Q_i)}^2. \quad (3.28)$$

By using estimate (3.28) and choosing the parameter $\theta \leq (C_q C_{inv,0,u})^{-1}$, we obtain

$$a_h(v_h, v_h) \geq \|v_h\|_h^2 - (\theta_i C_{inv,0,u} C_q)/2 \sum_{i=1}^N \|\nabla_x v_h\|_{L_2(Q_i)}^2$$

$$\geq (1 - (\theta_i C_{inv,0,u} C_q)/2)\|v_h\|_h^2 \geq \frac{1}{2}\|v_h\|_h^2.$$

$\square$

## 4    Discretization Error Analysis

In this section, we present *a priori* error estimates for the space-time scheme. We will require the following the following trace inequality.

**Lemma 5.** *Let $Q_i = \Phi_i(\widehat{Q})$ for $i = 1, \ldots, N$. Then the patch-wise scaled trace inequality*

$$\|v\|_{L_2(\partial Q_i)} \leq C_{t,u} h_i^{-1/2}\left(\|v\|_{L_2(Q_i)} + h_i|v|_{H^1(Q_i)}\right), \qquad (4.1)$$

*holds for all $v \in H^1(Q_i)$, where $h_i$ denotes the maximum mesh size in the physical domain, and $C_{t,u}$ is a positive constant that only depends on the shape regularity of the mapping $\Phi_i$.*

Next, we show the uniform boundedness of the discrete bilinear form $a_h(\cdot, \cdot)$ on $V_{0h,*} \times V_{0h}$, where the space $V_{0h,*} = H^2(Q) + V_{0h}$ is equipped with the norm

$$\|v\|_{h,*}^2 = \|v\|_h^2 + \sum_{i=1}^N (\theta_i h_i)^{-1}\|v\|_{L_2(Q_i)}^2 + (\theta_i h_i)^2\|\partial_t \nabla_x v\|_{L_2(Q_i)}^2$$

$$+ \sum_{F_{ij} \in \mathcal{F}_I} \|\{v\}^{up}\|_{L_2(F_{ij})}^2 + \theta_i h_i \sum_{F_{ij} \in \mathcal{F}_I} \|\{\nabla_x v\}^{up}\|_{L_2(F_{ij})}^2 + \sum_{i=1}^N h_i \|\nabla_x v_i\|_{L_2(\partial Q_i)}^2.$$

$$(4.2)$$

In order to prove boundedness of the bilinear form $a_h(\cdot, \cdot)$, we need two auxiliary inequalities presented in the following lemma.

**Lemma 6.** *For positive parameters $\delta_1$ and $\delta_2$, and for $F_{ij} \in \mathcal{F}_I$, the estimates*

$$\left| \int_{F_{ij}} \{\nabla_x u\}[\![v_h]\!]_x \, ds \right| \leq \left( \frac{\delta_1^{-1} h_i}{4}(\|\nabla_x u_i\|_{L_2(F_{ij})}^2 + \|\nabla_x u_j\|_{L_2(F_{ij})}^2) \right)^{\frac{1}{2}}$$

$$\times \left( \frac{\delta_1}{h_i} \| [\![ v_h ]\!]_x \|^2_{L_2(F_{ij})} \right)^{\frac{1}{2}}, \tag{4.3}$$

$$\left| \theta_i h_i \int_{F_{ij}} \{ \nabla_x u \} [\![ \partial_t v_h ]\!]_x \, ds \right| \leq \left( \frac{\delta_2^{-1} \theta_i h_i}{4} ( \| \nabla_x u_i \|^2_{L_2(F_{ij})} + \| \nabla_x u_j \|^2_{L_2(F_{ij})} ) \right)^{\frac{1}{2}}$$

$$\times \left( \delta_2 \theta_i h_i \| [\![ \partial_t v_h ]\!]_x \|^2_{L_2(F_{ij})} \right)^{\frac{1}{2}}, \tag{4.4}$$

hold for all $u \in V_{0h,*}$ and for all $v_h \in V_{0h}$, where $h_i$ is the maximum mesh size of patch $Q_i$.

*Proof.* The proof of (4.3) is obtained by using the Cauchy-Schwarz and triangle inequalities and the proof of (4.4) is obtained by means of the Cauchy-Schwarz inequality, See [21, Lemma 6.12]. $\square$

**Lemma 7.** *The discrete bilinear form $a_h(\cdot, \cdot)$, defined by (3.21), is uniformly bounded on $V_{0h,*} \times V_{0h}$, i.e., there exists a positive constant $\mu_b$ that does not depend on $h$ such that*

$$|a_h(u, v_h)| \leq \mu_b \| u \|_{h,*} \| v_h \|_h, \ \forall \, u \in V_{0h,*}, \forall \, v_h \in V_{0h}. \tag{4.5}$$

*Proof.* Since most of the terms are bounded by using Cauchy-Schwarz's inequality, we will only show the estimates for the other terms. For the next terms, we proceed by using Cauchy-Schwarz inequality as follows

$$\sum_{F_{ij} \in \mathcal{F}_I} \int_{F_{ij}} \{ \nabla_x v_h \} [\![ u ]\!]_x \, ds \leq \left( \sum_{F_{ij} \in \mathcal{F}_I} h_i \delta_1^{-1} \| \{ \nabla_x v_h \} \|^2_{L_2(F_{ij})} \right)^{\frac{1}{2}} \left( \sum_{F_{ij} \in \mathcal{F}_I} h_i^{-1} \delta_1 \| [\![ u ]\!]_x \|^2_{L_2(F_{ij})} \right)^{\frac{1}{2}}. \tag{4.6}$$

Considering the first term in (3.21), we apply the quasi-uniformity property (2.8) together with the inverse inequality (3.26) to obtain

$$\sum_{F_{ij} \in \mathcal{F}_I} h_i \delta_1^{-1} \| \{ \nabla_x v_h \} \|^2_{L_2(F_{ij})} \leq \sum_{i=1}^{N} (\delta_1^{-1} C_q C_{inv,0,u})/2 \| \nabla_x v_{h,i} \|^2_{L_2(Q_i)}. \tag{4.7}$$

Using the techniques of the proof of (4.6), we estimate the next term as follows:

$$\theta_i h_i \sum_{F_{ij} \in \mathcal{F}_I} \int_{F_{ij}} \{ \nabla_x v_h \} [\![ \partial_t u ]\!]_x \, ds$$

$$\leq \left( \sum_{F_{ij} \in \mathcal{F}_I} \delta_2^{-1} \theta_i h_i \| \{ \nabla_x v_h \} \|_{L_2(F_{ij})}^2 \right)^{\frac{1}{2}} \left( \sum_{F_{ij} \in \mathcal{F}_I} \theta_i h_i \delta_2 \| [\![ \partial_t u ]\!]_x \|_{L_2(F_{ij})}^2 \right)^{\frac{1}{2}}.$$

The first term is estimated similarly to (4.7) as follows

$$\sum_{F_{ij} \in \mathcal{F}_I} \delta_2^{-1} \theta_i h_i \| \{ \nabla_x v_h \} \|_{L_2(F_{ij})}^2 \leq \sum_{i=1}^N (\theta_i \delta_2^{-1} C_q C_{inv,0,u}) / 2 \| \nabla_x v_{h,i} \|_{L_2(Q_i)}^2.$$

Using Lemma 6, the next terms are estimated by means of the quasi-uniformity property (2.8) as follows

$$\sum_{F_{ij} \in \mathcal{F}_I} \int_{F_{ij}} \{ \nabla_x u \} [\![ v_h ]\!]_x \, ds$$

$$\leq \left( \delta_1^{-1} / 4 \sum_{i=1}^N C_q h_i \| \nabla_x u_i \|_{L_2(\partial Q_i)}^2 \right)^{1/2} \left( \sum_{F_{ij} \in \mathcal{F}_I} \frac{\delta_1}{h_i} \| [\![ v_h ]\!]_x \|_{L_2(F_{ij})}^2 \right)^{1/2},$$

$$\theta_i h_i \sum_{F_{ij} \in \mathcal{F}_I} \int_{F_{ij}} \{ \nabla_x u \} [\![ \partial_t v_h ]\!]_x \, ds$$

$$\leq \left( \delta_2^{-1} \theta_i / 4 \sum_{i=1}^N C_q h_i \| \nabla_x u_i \|_{L_2(\partial Q_i)}^2 \right)^{1/2} \left( \sum_{F_{ij} \in \mathcal{F}_I} \delta_2 \theta_i h_i \| [\![ \partial_t v_h ]\!]_x \|_{L_2(F_{ij})}^2 \right)^{1/2}.$$

Combining the terms from above and using Cauchy Schwarz's inequality, we obtain

$$|a_h(u, v_h)| \leq \mu_b \| u \|_{h,*} \| v_h \|_h,$$

where $\mu_b = 2 \sqrt{\max \{ C_{inv,0,u} C_q (\delta_1^{-1} + \theta_i \delta_2^{-1}) / 4 + 1, (\delta_1^{-1} + \theta_i \delta_2^{-1}) / 8 \}}$.    $\square$

We recall the local interpolation error estimate as presented in [4], then, using a patch-wise interpolation error estimate, we provide an interpolation error estimate in the energy norm.

**Proposition 1.** *Given the integers $l_i$ and $s_i$ such that $0 \leq l_i \leq s_i \leq p_i + 1$, for a function $v \in H^{s_i}(Q_i)$, then*

$$\sum_{K \in \mathcal{K}_h} |v - \Pi_h v|_{H^{l_i}(K)} \leq C_s h_i^{s_i - l_i} \| v \|_{H^{s_i}(Q_i)}, \tag{4.8}$$

*where $h_i$ denotes the maximum mesh-size parameter in the physical domain and the constant $C_s$ only depends on $l_i, s_i$ and $p_i$, the shape regularity of the physical domain $Q_i$ described by the mapping $\Phi_i$ and, in particular, $\nabla \Phi_i$.*

*Proof.* See, e.g. [26, Proposition 3.1].

If the multiplicity of the inner knots is not larger than $p_i + 1 - l_i$, then $\Pi_h v \in V_h \cap H^{l_i}(Q_i)$, and estimate (4.8) yields a global estimate.

**Proposition 2.** *Let us assume that the multiplicity of the inner knots is not larger than $p_i + 1 - l_i$. Given the integers $l$ and $s$ such that $0 \leq l_i \leq s_i \leq p_i + 1$, there exist a positive constant $C_s$ such that for a function $v \in H^{s_i}(Q_i)$*

$$|v - \Pi_h v|_{H^{l_i}(Q_i)} \leq C_s h^{(s_i - l_i)} \|v\|_{H^{s_i}(Q)}, \qquad (4.9)$$

*where $h_i$ denotes the maximum mesh-size parameter in the physical domain and the generic constant $C_s$ only depends on $l_i, s_i$ and $p_i$, the shape regularity of the physical domain $Q_i$ described by the mapping $\Phi_i$ and, in particular, $\nabla \Phi_i$.*

*Proof.* See e.g. [26, Proposition 3.2].

**Lemma 8.** *Let $\mathbf{s}$ be a positive integer with $\mathbf{s} \geq 2$ and $v \in H^{\mathbf{s}}(Q, \mathcal{T}_h)$. Then, there exists a projective operator $\Pi_h$ such that $\Pi_h v \in V_{0h}$ and generic positive constants $C_0, C_1, C_2, C_3, C_4$ and $C_5$ such that*

$$\|\nabla(v - \Pi_{h,i} v)\|^2_{L_2(\partial Q_i)} \leq C_0 h_i^{2(r_i - 1) - 1} \|v\|^2_{H^{r_i}(Q_i)}, \qquad (4.10)$$

$$\sum_{F_{ij} \in \mathcal{F}_I} h \| [\![ \nabla_x (v - \Pi_h v) ]\!]_t \|^2_{L_2(F_{ij})} \leq C_1 \sum_{i=1}^N h_i^{2(r_i - 1)} \|v\|^2_{H^{r_i}(Q_i)}, \qquad (4.11)$$

$$\sum_{F_{ij} \in \mathcal{F}_I} \delta_2 \theta_i h_i \| [\![ \partial_t (v - \Pi_h v) ]\!]_x \|^2_{L_2(F_{ij})} \leq C_2 \sum_{i=1}^N h_i^{2(r_i - 1)} \|v\|^2_{H^{r_i}(Q_i)}, \qquad (4.12)$$

$$\sum_{F_{ij} \in \mathcal{F}_I} \frac{\delta_1}{h} \| [\![ v - \Pi_h v ]\!]_x \|^2_{L_2(F_{ij})} \leq C_3 \sum_{i=1}^N h_i^{2(r_i - 1)} \|v\|^2_{H^{r_i}(Q_i)}, \qquad (4.13)$$

$$\sum_{F_{ij} \in \mathcal{F}_I} \| \{ v - \Pi_h v \}^{up} \|^2_{L_2(F_{ij})} \leq C_4 \sum_{i=1}^N h_i^{2r_i - 1} \|v\|^2_{H^{r_i}(Q_i)}, \qquad (4.14)$$

$$\sum_{F_{ij} \in \mathcal{F}_I} \theta_i h_i \| \{ \nabla_x (v - \Pi_h v) \}^{up} \|^2_{L_2(F_{ij})} \leq C_5 \sum_{i=1}^N h_i^{2(r_i - 1)} \|v\|^2_{H^{r_i}(Q_i)}, \qquad (4.15)$$

*where $r_i = \min\{s_i, p_i + 1\}$, and $h_i$ is the mesh size of the patch $Q_i$ in the physical domain, $p_i$ is polynomial degree of the NURBS, $\delta_1$ and $\delta_2$ are positive constants and the generic constants are independent of the mesh size.*

*Proof.* By using the trace inequality of Lemma 5 and the approximation estimate of Proposition 2, we proceed with the proof of (4.10) as follows

$$\|\nabla(v - \Pi_{h,i}v)\|^2_{L_2(\partial Q_i)}$$
$$\leq C^2_{t,u}h^{-1}_i\left(\|\nabla(v - \Pi_{h,i}v)\|^2_{L_2(Q_i)} + h^2_i|\nabla(v - \Pi_{h,i}v)|^2_{H^1(Q_i)}\right)$$
$$\leq C_s C^2_{t,u}h^{-1}_i\left(h^{2(r_i-1)}_i\|v\|^2_{H^{r_i}(Q_i)} + h^2_i h^{2(r-2)}_i\|v\|^2_{H^{r_i}(Q_i)}\right)$$
$$\leq C_0 h^{2(r_i-1)-1}_i\|v\|^2_{H^{r_i}(Q_i)}, \quad C_0 = 2C_s C^2_{t,u}.$$

Since $\nabla = (\nabla_x, \partial_t)^\top$, the approximation (4.10) also holds for both components of the derivative. We will particularly need the following

$$\|\nabla_x(v - \Pi_{h,i}v)\|^2_{L_2(\partial Q_i)} \leq C_0 h^{2(r_i-1)-1}_i\|v\|^2_{H^{r_i}(Q_i)},$$
$$\|\partial_t(v - \Pi_{h,i}v)\|^2_{L_2(\partial Q_i)} \leq C_0 h^{2(r_i-1)-1}_i\|v\|^2_{H^{r_i}(Q_i)}. \tag{4.16}$$

Next, by using Definition 1, (4.16) and (4.10) together with (2.8), the proofs of (4.11)–(4.13) are obtained. Since the last two estimates follow the same proof technique, we will only show the proof of (4.15). We estimate the upwind term (4.14). By using Definition 1 and (4.10) we proceed as follows

$$\sum_{F_{ij}\in\mathcal{F}_I}\|\{v - \Pi_h v\}^{up}\|^2_{L_2(F_{ij})}$$
$$\leq 2\sum_{i=1}^N C^2_{t,u}h^{-1}_i\left(\|v - \Pi_{h,i}v\|^2_{L_2(Q_i)} + h^2_i|v - \Pi_{h,i}v|^2_{H^1(Q_i)}\right)$$
$$\leq 2C_s C^2_{t,u}\sum_{i=1}^N h^{-1}_i\left(h^{2r}_i\|v\|^2_{H^{r_i}(Q_i)} + h^2_i h^{2(r_i-1)}_i\|v\|^2_{H^{r_i}(Q_i)}\right)$$
$$\leq 4C_s C^2_{t,u}\sum_{i=1}^N h^{2r_i-1}_i\|v\|^2_{H^{r_i}(Q_i)} \leq C_4 h^{2r_i-1}\|v\|^2_{H^{r_i}(Q_i)}.$$

$$\square$$

Next, we obtain estimates of the approximation error with respect to the discrete norms $\|\cdot\|_h$ and $\|\cdot\|_{h,*}$.

**Lemma 9.** *Let $s_i, i = 1,\ldots,N$ be a positive integer with $s_i \geq 2$ and $v \in H^{s_i}(Q_i), i = 1,\ldots,N$. Then there exist a projection $\Pi_h v \in V_{0h}$ and generic positive constants $C_6$ and $C_7$ such that*

$$\|v - \Pi_h v\|^2_h \leq C_6\sum_{i=1}^N h^{2(r_i-1)}_i\|v\|^2_{H^{r_i}(Q_i)}, \tag{4.17}$$

$$\|v - \Pi_h v\|_{h,*}^2 \leq C_7 \sum_{i=1}^{N} h^{2(r_i-1)} \|v\|_{H^{r_i}(Q_i)}^2, \qquad (4.18)$$

where $r_i = \min\{s_i, p_i + 1\}$ and $h_i$ is the mesh-size of $Q_i$, $p_i$ is polynomial degree of the NURBS and the generic constants $C_6$ and $C_7$ are independent of the mesh size $h_i$.

*Proof.* Using the definition of the norms $\|\cdot\|_h$ and $\|\cdot\|_{h,*}$ together with Lemma 8, we are able to obtain the proof of the statement.  □

Next, we show the relationship between our model problem (2.1) respectively (2.2) and the space-time dG-IgA variational identity (3.20) is given by the consistency theorem.

**Theorem 1.** *If the solution of the variational problem (2.2) belongs to $H^2(Q)$ then the solution $u$ satisfies the space-time dG-IgA consistency identity (3.20).*

*Proof.* Following the same arguments as in [20, Lemma 12], we obtain that

$$\partial_t u - \Delta u = f \quad \text{in } L_2(Q) \qquad \text{and} \qquad u = 0 \quad \text{in } L_2(\Sigma_0). \qquad (4.19)$$

We multiply the differential equation of (4.19) with the test function $v_h + \theta h \partial_t v_h$ for $v_h \in V_{0h}$ and integrate over the space-time domain $Q$. Since $u \in H^2(Q)$ and $p \geq 2$, we can apply all the derivations as we did in Section 3.  □

Finally, we present the main result for this chapter, namely the *a priori* error discretization error estimate in the discrete norm $\|\cdot\|_h$.

**Theorem 2.** *Let $u$ be the exact solution of our model problem (2.2), and let $u_h$ be the solution to the space-time dG-IgA scheme (3.20). Then the discretization error estimate*

$$\|u - u_h\|_h \leq C \sum_{i=1}^{N} h_i^{r_i-1} \|v\|_{H^{r_i}(Q_i)}, \qquad (4.20)$$

*holds, where $r_i = \min\{s_i, p_i + 1\}$, $C$ is a generic positive constant independent of the mesh size $h_i$ and $p_i$ denotes the underlying polynomial degree of the NURBS patch $Q_i$.*

*Proof.* Subtracting the space-time dG-IgA scheme

$$a_h(u_h, v_h) = \ell_h(v_h), \quad \forall v_h \in V_{0h},$$

from the consistency identity

$$a_h(u, v_h) = \ell_h(v_h), \quad \forall v_h \in V_{0h},$$

we obtain the so-called Galerkin orthogonality

$$a_h(u - u_h, v_h) = 0, \quad \forall v_h \in V_{0h}, \tag{4.21}$$

that is crucial for the discretization error estimate. Using now the triangle inequality, we can estimate the discretization error $u - u_h$ as follows

$$\|u - u_h\|_h \leq \|u - \Pi_h u\|_h + \|\Pi_h u - u_h\|_h. \tag{4.22}$$

The first term is estimated by Lemma 9. The estimation of the second term on the right-hand side of (4.22) follows by using the fact that $\Pi_h u - u_h \in V_{0h}$, the ellipticity of the bilinear form $a_h(\cdot, \cdot)$ as was shown in Lemma 4, the Galerkin orthogonality (4.21), and the boundedness of the discrete bilinear form Lemma 7, we can derive the following estimates

$$\begin{aligned}
\mu_c \|\Pi_h u - u_h\|_h^2 &\leq a_h(\Pi_h u - u_h, \Pi_h u - u_h) = a_h(\Pi_h u - u, \Pi_h u - u_h) \\
&\leq \mu_b \|\Pi_h u - u\|_{h,*} \|\Pi_h u - u_h\|_h.
\end{aligned}$$

Hence, we have

$$\|\Pi_h u - u_h\|_h \leq (\mu_b/\mu_c) \|\Pi_h u - u\|_{h,*}. \tag{4.23}$$

Inserting (4.23) into the triangle inequality (4.22) and using the estimates (4.18) and (4.17) from Lemma 9, we have

$$\begin{aligned}
\|u - u_h\|_h &\leq \|u - \Pi_h u\|_h + \|\Pi_h u - u_h\|_h \\
&\leq \|u - \Pi_h u\|_h + (\mu_b/\mu_c) \|\Pi_h u - u\|_{h,*} \\
&\leq (C_6 + (\mu_b/\mu_c)C_7) \sum_{i=1}^{N} h_i^{r_i - 1} \|u\|_{H^{r_i}(Q_i)},
\end{aligned}$$

which proves the discretization error estimate (4.20) with $C = (C_6 + (\mu_b/\mu_c)C_7)$. $\quad\square$

## 5  Numerical Results

The numerical results presented below have been performed in G+SMO, see [14]. We used the sparse direct solver SuperLU to solve the resulting linear system of space-time dG-IgA equations. We present numerical results for two and three dimensional spatial computational domains with the parameter $\theta_i = 0.1, i = 1, \ldots, N$ in all our numerical experiments in the space-time $dG-$norm as well as $L_2$-norm. We assume in all our experiments that the B-spline degree in both space-direction and time-direction are equal. We also note that the B-spline degree on all the patches are equal and denoted by $p$ i.e. $p = p_i, i = 1, \ldots, N$. The penalty parameters are chosen as $\delta_1 = \delta_2 = 2(p + d + 1)(p + 1)$. In our numerical experiments on the convergence behavior of the space-time dG-IgA methods proposed, we have used B-splines of polynomial degrees $p = 1, 2, 3, 4$ and computed the rate of convergence of successive mesh refinement by means of the formula $\log_2(\| \cdot \|_{h,i+1} / \| \cdot \|_{h,i})$.

### 5.1  Two dimensional space-time computational domain

We consider a space-time computational domain $Q = \{(x, t) \in \mathbb{R}^2 : x \in \Omega(t), t \in (0, 2)\}$ with a one-dimensional moving spatial computational domain of the form $\Omega(t) = \{x \in \mathbb{R} : a(t) < x < b(t)\}, t \in (0, 2)$, where $a(t) = -t(1 - t)/2$ and $b(t) = 1 - t(1 - t)/2$ for $t \in [0, 1]$, whereas $a(t) = -(t - 1)(t - 2)/2$ and $b(t) = (t^2 - 3t + 4)/2$ for $t \in [1, 2]$, leading to the decomposition of the space-time cylinder $\overline{Q} = \overline{Q}_2 \cup \overline{Q}_1$ into the two non-overlapping patches $Q_2 = \{(x, t) \in \mathbb{R}^2 : x \in \Omega(t), t \in (0, 1)\}$ and $Q_1 = \{(x, t) \in \mathbb{R}^2 : x \in \Omega(t), t \in (1, 2)\}$, respectively. The two patches $Q_1$ and $Q_2$ can also be represented by knot vectors such that $Q_1 = \{\Xi_1 := \{0, 0, 1, 1\}, \Xi_2 := \{0, 0, 0, 1, 1, 1\}\}$ and $Q_2 = \{\Xi_1 := \{0, 0, 1, 1\}, \Xi_2 := \{0, 0, 0, 1, 1, 1\}\}$ with the corresponding control points $\mathbf{P}_{1,1}^1 = (0, 1)$, $\mathbf{P}_{2,1}^1 = (1, 1)$, $\mathbf{P}_{2,2}^1 = (0.25, 1.5)$, $\mathbf{P}_{2,3}^1 = (0.75, 1.5)$, $\mathbf{P}_{1,3}^1 = (0, 2)$ and $\mathbf{P}_{1,2}^1 = (1, 2)$ and $\mathbf{P}_{1,1}^2 = (0, 0)$, $\mathbf{P}_{2,1}^2 = (1, 0)$, $\mathbf{P}_{2,2}^2 = (-0.25, 0.5)$, $\mathbf{P}_{2,3}^2 = (1.25, 0.5)$, $\mathbf{P}_{1,3}^2 = (0, 1)$ and $\mathbf{P}_{1,2}^2 = (1, 1)$, respectively. We again solve our model problem (2.1), and choose the data such that the solution is given by $u(x, t) = \sin(\pi x) \sin(\pi t)$, i.e. $f(x, t) = \partial_t u(x, t) - \Delta u(x, t) = (\pi \sin(\pi x))(\cos(\pi t) + \pi \sin(\pi t))$ in $Q$, $u_0 = 0$ on $\overline{\Omega}$, and $u(x, t) = \sin(\pi x) \sin(\pi t)$ on $\Sigma$. Thus, the compatibility condition between boundary and initial conditions holds. The convergence behavior of the space-time dG-IgA scheme with respect to the discrete norm $\| \cdot \|_h$ is shown in Figure 4 by a series of $h$-refinement and by using

B-splines of polynomial degrees $p = 1, 2, 3, 4$. After some saturation, we observe the optimal convergence rate $O(h^p)$ for $p \geq 2$ as theoretically predicted by Theorem 2 for smooth solutions. Moreover, Figure 4 (right) shows the $L_2$ errors and the corresponding rates for the same setting. We see that the $L_2$ rates are asymptotically optimal for $p \geq 2$ as well, i.e. they behave like $O(h^{p+1})$. For $p = 1$, we also observe the optimal rate in the dG-norm $\| \cdot \|_h$, whereas the $L_2$-rate does not reach the optimal order 2.

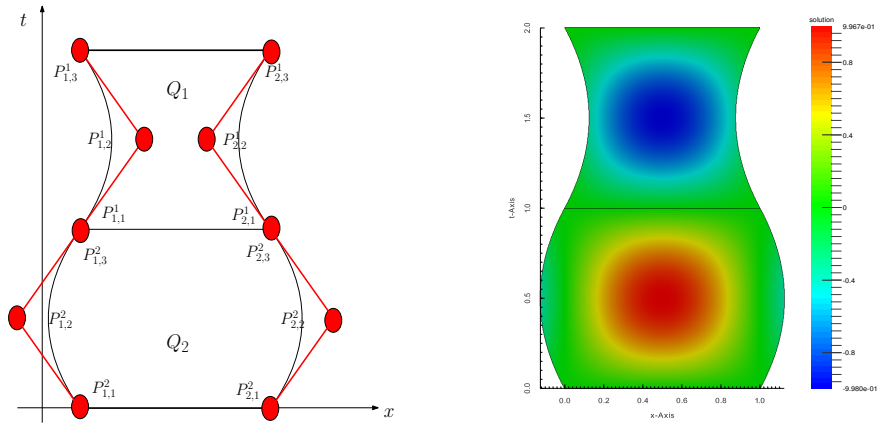

**Fig. 3.** Multipatch space-time computational domain with control points (left) and solution contours (right) for Example 5.1.

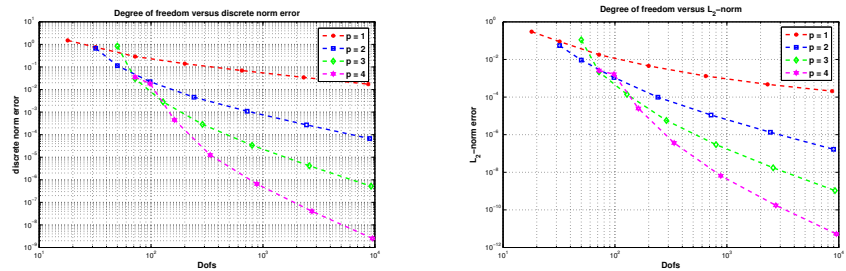

**Fig. 4.** Moving spatial domain $Q \subset \mathbb{R}^2$ : log-log plot of $\| \cdot \|_h$ (left) and $\| \cdot \|_{L_2(Q)}$ (right) versus Degree of freedom (Dofs).

## 5.2   Three dimensional space-time computational domain.

We again consider a three-dimensional space-time computational domain $Q = \{(x,t) \in \mathbb{R}^3 : x \in \Omega(t), t \in (0,2)\} \subset \mathbb{R}^3$. The space-time cylinder is decomposed as follows $\overline{Q} = \overline{Q}_2 \cup \overline{Q}_1$ into the two non-overlapping patches $Q_2 = \{(x,t) \in \mathbb{R}^3 : x \in \Omega(t)\}, t \in (0,1)$ and $Q_1 = \{(x,t) \in \mathbb{R}^3 : x \in \Omega(t)\}, t \in (1,2)$. The two patches $Q_1$ and '$Q_2$ can also be represented by knot vectors such that $Q_1 = \{\Xi_1 := \{0,0,1,1\}, \Xi_2 := \{0,0,1,1\}, \Xi_3 := \{0,0,0,1,1,1\}\}$ and $Q_2 = \{\Xi_1 := \{0,0,1,1\}, \Xi_2 := \{0,0,1,1\}, \Xi_3 := \{0,0,0,1,1,1\}\}$ with the corresponding control points $\mathbf{P}^1_{i_1,i_2,i_3}$ and $\mathbf{P}^2_{i_1,i_2,i_3}$, presented in Table 1. We solve our model problem (2.1), and again choose the data such that the solution is given by $u(x,t) = \sin(\pi x)\sin(\pi t)$, i.e. $f(x,t) = \partial_t u(x,t) - \Delta u(x,t) = (\pi \sin(\pi x))(\cos(\pi t) + \pi \sin(\pi t))$ in $Q$, $u_0 = 0$ on $\overline{\Omega}$, and $u(x,t) = \sin(\pi x)\sin(\pi t)$ on $\Sigma$. Thus, the compatibility condition between boundary and initial conditions holds. The convergence behavior of the space-time dG-IgA scheme with respect to the space-time dG norm $\|\cdot\|_h$ is shown in Figure 7 (left) by a series of $h$-refinement and by using B-splines of polynomial degrees $p = 1,2,3,4$. After some saturation, we observe the optimal convergence rate $O(h^p)$ for $p \geq 2$ as theoretically predicted by Theorem 2 for smooth solutions. Moreover, Figure 7 (right) shows the $L_2$ errors and the corresponding rates for the same setting. We see that the $L_2$ rates are asymptotically optimal for $p \geq 2$ as well, i.e. they behave like $O(h^{p+1})$. For $p = 1$, we also observe the optimal rate in the discrete norm, whereas the $L_2$-rate does not reach the optimal order 2.

| $i_1$ | $i_2$ | $i_3$ | $\mathbf{P}^1_{i_1,i_2,i_3}$ | $\mathbf{P}^2_{i_1,i_2,i_3}$ |
|---|---|---|---|---|
| 1 | 1 | 1 | $(0,0,1)$ | $(0,0,0)$ |
| 1 | 1 | 2 | $(0.25,0,1.5)$ | $(-0.25,0,0.5)$ |
| 1 | 1 | 3 | $(0,0,2)$ | $(0,0,1)$ |
| 1 | 2 | 1 | $(0,1,1)$ | $(0,1,0)$ |
| 1 | 2 | 2 | $(0.25,1,1.5)$ | $(-0.25,1,0.5)$ |
| 1 | 2 | 3 | $(0,1,2)$ | $(0,1,1)$ |
| 2 | 1 | 1 | $(1,0,1)$ | $(1,0,0)$ |
| 2 | 1 | 2 | $(0.75,0,1.5)$ | $(1.25,0,0.5)$ |
| 2 | 1 | 3 | $(1,0,2)$ | $(1,0,1)$ |
| 2 | 2 | 1 | $(1,1,1)$ | $(1,1,0)$ |
| 2 | 2 | 2 | $(0.75,1,1.5)$ | $(1.25,1,0.5)$ |
| 2 | 2 | 3 | $(1,1,2)$ | $(1,1,1)$ |

**Table 1.** The control points $\mathbf{P}^1_{i_1,i_2,i_3}$ and $\mathbf{P}^2_{i_1,i_2,i_3}$ corresponding to patches $Q_1$ and $Q_2$ as depicted in Figure 6.

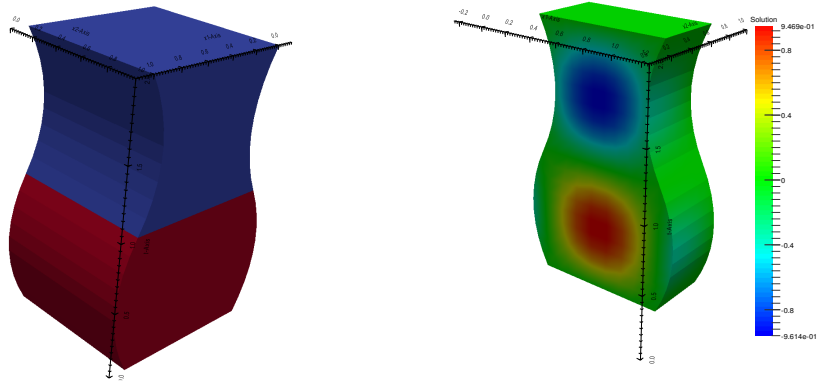

**Fig. 5.** The multi-patch space-time computational domain (left) and the solution contours for Example 5.2 (right).

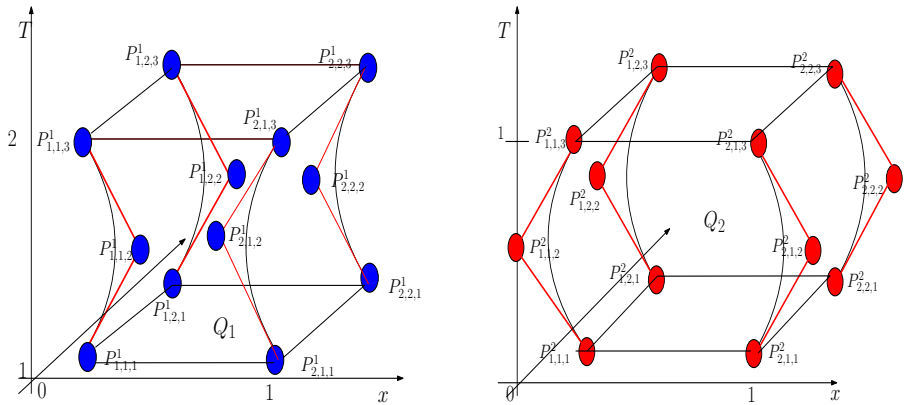

**Fig. 6.** The control points of the decomposed space-time computational domain consisting of patches $Q_1$ (left) and $Q_2$ (right).

## 6   Conclusion

In this paper, we presented a space-time multi-patch discontinuous Galerkin isogeometric analysis (dG-IgA) for parabolic evolution problems. Parabolic evolution problems are usually solved by time-stepping scheme. However, by using a time-upwind test function, we derive a stable or elliptic space-time multi-patch scheme. We applied discontinuous Galerkin methodology in both space and time and derived a stable space-time dG scheme. We showed optimal *a priori* error estimates with respect to a space-time dG norm $\| \cdot \|_h$. Finally, we presented convincing numerical experiments for the scheme including two and three dimensional computational domains. The current scheme can be combined with surface PDEs as presented in [19]. An extension of the analysis to include low regularity space-time solutions, i.e., $u \in W^{2,q}(Q)$ following ideas presented in the monograph [9] will be considered in the future.

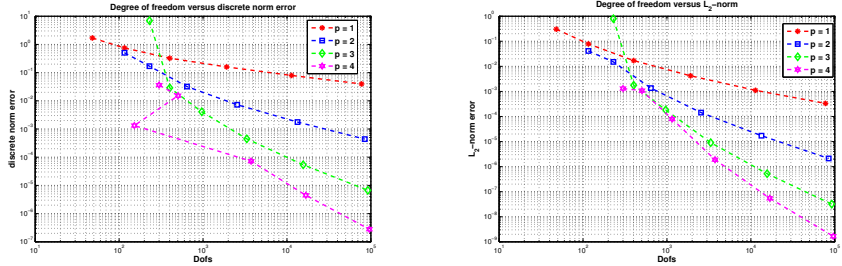

**Fig. 7.** Moving spatial domain $Q \subset \mathbb{R}^3$ : log-log plot of $\| \cdot \|_h$ (left) and $\| \cdot \|_{L_2(Q)}$ (right) versus Degree of freedom (Dofs).

## Acknowledgement

The author would like to thank Martin Neumüller and Ulrich Langer for several discussions during the preparation of the article.